\documentclass[a4paper]{article}
\usepackage[toc,page]{appendix}

\pdfoutput=1

\usepackage{titlesec}
\usepackage[latin1]{inputenc}
\usepackage{enumerate}
\usepackage{graphicx}
\usepackage{subfigure}
\usepackage{rotating}
\usepackage{amsmath,amsthm,amsfonts,amssymb,color,graphicx}
\usepackage[usenames,dvipsnames]{pstricks}
\usepackage{epsfig}
\usepackage{pst-grad}
\usepackage{pst-plot}
\usepackage[final]{pdfpages}
\usepackage{seqsplit}
\usepackage{physics}
\usepackage{mathtools}
\usepackage{listings}

\usepackage{comment}
\usepackage{chngcntr}\counterwithin{equation}{section}

\usepackage{braket}
\usepackage[breaklinks]{hyperref}

\newtheorem{theorem}{Theorem}[section]
\newtheorem{proposition}[theorem]{Proposition}
\newtheorem{lemma}[theorem]{Lemma}

\newtheorem{corollary}[theorem]{Corollary}

\def\Z{\mathbb{Z}}
\def\N{\mathbb{N}}
\def\Q{\mathbb{Q}}
\def\R{\mathbb{R}}

\def\D{\mathbb{D}}

\newcommand{\leg}[2]{({\scriptstyle{\frac{#1}{#2}}})}

\def\lra{{\longrightarrow}}

\newcommand{\cheb}[2]{T_{#1}\left({#2}\right)}
\newcommand{\chsh}[2]{T^*_{#1}\left({#2}\right)}
\newcommand{\chebh}[1]{T_{#1}}
\newcommand{\chshh}[1]{T^*_{#1}}

\newtheorem*{theorem*}{Theorem}

\def\ort{\textrm{ord}}
\def\qD{\sqrt{D}}
\def\qrD{ \Q( \qD ) }
\def\fp{\mathfrak p}

\def\Gal{\mathrm Gal}
\def\cR{\mathbf{R}}
\def\cJ{\mathcal{J}}

\def\TK{\mathcal{T}_K}

\def\blah{$\left(\mathbf{W}\right)$}

\begin{document}
\author{Gary McConnell}

\title{Some new infinite families of non-$p$-rational real quadratic fields}
\maketitle

\begin{abstract}
Fix a finite collection of primes~$\{ p_j \}$, not containing $2$ or $3$. 
Using some observations which arose from attempts to solve the SIC-POVMs problem in quantum information, we give a simple methodology for constructing an infinite family of simultaneously non-$p_j$-rational real quadratic fields, unramified above any of the $p_j$. 
Alternatively these may be described as infinite sequences of instances of~$\qrD$, for varying~$D$, where every $p_j$ is a $k$-Wall-Sun-Sun prime, or equivalently a generalised Fibonacci-Wieferich prime. 

One feature of these techniques is that they may be used to yield fields $K= \qrD$ for which a $p$-power cyclic component of the torsion group of the Galois groups of the maximal abelian pro-$p$ extension of $ \qrD$ unramified outside primes above $p$, is of size $p^a$ for $a\geq1$ arbitrarily large. 
\end{abstract}

\vskip 5mm

{\raggedleft \bf Dedicated to the memory of John Coates \rm\par}
\vskip 5mm

\section{Introduction}
In \cite{AFMY} some relations were conjectured between problems in quantum state space geometry in complex Hilbert space of dimension $d$ -- amounting to finding maximal sets of equiangular lines -- and ray class fields of finite modulus $d$ of real quadratic fields $K= \qrD$; where $D$ is the square-free part of $d^2-2d-3$. 
Subsequently the links between these disparate objects have been further strengthened by the discovery that square roots of Galois conjugates of Stark units associated to these very ray class fields actually seem to provide the `complex angles' which generate the geometry itself~\cite{kopp, ABGHM, BGM}. 

This peculiar construct means that attached to any particular real quadratic field $K$ is an \emph{infinite tower of dimensions} $d_\ell = d_\ell(D)$ which is monotonically increasing in $\ell$. 
It is furthermore partially ordered by divisibility among the $\ell$, subject to a shift to an entirely new set of primes --- with the possible exception of the prime $3$ itself --- above every successive power~$\ell=3^s$ of $3$. 
We give an outline of the main ingredients of this system in the appendix to this paper, but it is introduced in \cite[\S3]{AFMY} and explained more fully in appendix~A of~\cite{BGM}. 

One offshoot of this work has been a slightly obtuse entry-point into some 
problems in elementary number theory 
regarding so-called $k$-Wall-Sun-Sun primes. 
In the course of conversations with John Coates in 2017 on these mysterious links, we 
observed that what might be called a second-order 
aspect of Leopoldt's `conjecture'~\cite{wash} for real 
quadratic fields, cropped up obliquely inside these dimension towers. 

Since the sequences which gave rise to these results appear not to have been studied before --- I thank Jeff Lagarias for pointing this out to me --- and because they give a new, simply-described infinite series of non-$p$-rational real quadratic fields for any specified prime $p$, it seems useful to have them available for reference. 
There will predictably be some considerable overlap with the infinite sequence of such fields in~\cite[\S6]{grashal} and~\cite[\S5.3]{grasBS}. 
The techniques we use are simpler to state; but they are concomitantly less general. 
Indeed, viewed from the perspective of a single fixed real quadratic field, we are restricted to examples from a subset of the primes which has proportion $\frac{3}{8}$ of all primes~\cite[\S4]{kopp}, which happens to have an internal structure enabling more control over some of the $p$-adic invariants. 
Presumably there is a similarly concise and easily-controlled criterion for the remaining $\frac{5}{8}$ of the primes, but we have yet to find it. 

\section{$p$-rationality for real quadratic fields}

\subsection{Definition of $p$-rationality}
\noindent
The class field-theoretic notion of the \emph{$p$-rationality} of a number field~$K$ was introduced in the context of so-called abelian $p$-ramification theory by Mohavvedi \cite{Mova} also with Nguyen Quang-Do \cite[\S2]{MQD} as a measure of the `normalcy' of a prime $p$ within a number field, assuming that Leopoldt's conjecture is known for that field, which of course for abelian extensions of $\Q$ it is \cite[theorem~5.25]{wash}. 
The idea is to measure the size of a finite Galois group $\TK$, the torsion subgroup of the maximal abelian pro-$p$-extension of $K$ unramified other than at primes above $p$, which is an invariant first defined as such by Gras in \cite{grascrelle1, grascrelle2, grasK2}, and which has three components \cite[\S2.3]{grasBS} in the case of totally real fields $K$. 
See also~\cite{jauthes}; this whole subject largely goes back to Coates' analytic formula in \cite[theorem 1.13]{coates}. 

Since we stick exclusively here to real quadratic fields, we describe the factors which can lead to a non-trivial group $\TK$ for a given $p$ and $K$. 
First there is  a component --- which disappears for $p\geq5$ --- to do with $p$-th roots of unity in the $p$-adic completion(s) of the field. 
Secondly we need the $p$-primary component of the ideal class group of $K$; with the one exception when $K=\Q(\sqrt{2})$ and $p=2$, which again does not concern us in this paper. 
Finally comes the component whose size is a measure of $p$-divisibility of global units inside the local completions at the primes above $p$ of $K$, which is the focus of this paper. 
Up to a $p$-adic unit, its order is also described as the normalised $p$-adic regulator \cite{grasNpR}, which in our case where $K$ is always a real quadratic field amounts to $\frac{1}{p}$ times the $p$-adic logarithm of a fundamental unit of $K$. 

If every one of these three components has order $1$ then $K$ is said to be $p$-rational; otherwise it is non-$p$-rational.

\subsection{New infinite sets of non-$p$-rational real quadratic fields}
The purpose of this note is to prove the following two theorems. 
We first give a concise version, which is easy to visualise; after that we state a more technical result with broader scope. 
Define a `discriminant' map $\D$ from $\N \setminus \{0,1,2, 3\}$ onto the set of square-free positive integers, by sending $N\in\N$ to the square-free part of $(N+1)(N-3)\ \ [\ = (N-1)^2 - 4 \ ]$. 
We shall explore this map $\D$ more fully in section \ref{scoo}; but suffice it to say that up to the factor of $4$ which we ignore, it takes our dimensions $d_\ell(D)$ back to the discriminant of the real quadratic field $\qrD$. 

\begin{theorem}\label{peequeue}
Let $p,q$ be distinct prime numbers, with $p\geq5$.
Choose integers $r\geq2$ and $t\geq0$ such that $q^t < \frac{2^p p^{p-r}}{3^p}$. 
Then the real quadratic field $ \Q( \sqrt{\D( p^rq^t )} )$ is not $p$-rational. 
\end{theorem}
\noindent
We shall prove this in section \ref{proof1}. 
An immediate corollary of the proof techniques is then the next result, proven in section \ref{proof2}.  

\def\bfr{{\mathbf{r}}}
\begin{theorem}\label{pqr}
For some $n\geq1$ choose $n$ distinct primes $p_j \geq 5$, each 
together with an integer $r_j \geq n+1$, defining an index vector~$\bfr = (r_1,r_2,\ldots,r_n)$. 
Call the set of all such index vectors $\mathbf{I}_n \subset \N^{n}$. 
Write $d^{(\bfr)} = \scriptstyle{ \prod\limits_{j=1}^n p_j^{r_j} }$ for the corresponding dimension, with $K^{(\bfr)}$ being the real quadratic field 
$ \Q\left( \sqrt{ \D \left( d^{(\bfr)} \right) }  \right) $. 
Then the set~\hbox{$\{ K^{(\bfr)} \colon \bfr \in \mathbf{I}_n\}$} 
is an infinite collection of distinct real quadratic fields all of which are non-$p_j$-rational for \rm every \it $j$. 
\end{theorem}

It is worth remarking that the same ideas, applied this time rather as \emph{restrictions} on the growth of $p$-powers, can give us infinite sequences of {candidates} for \emph{$p$-rational} --- as opposed to non-$p$-rational --- real quadratic fields, contingent on being able to control the $p$-part of the ideal class group. 

The general expectation among the cognoscenti seems to be that `almost all' fields are $p$-rational at `almost all' primes. 
But not much is known beyond the results of Silverman \emph{et al} \cite{silverman,khare,Bou1,maire1,maire2} which prove that there are infinitely many $p$-rational fields for any given $p$ within certain classes of totally real fields, under the hypothesis of the $abc$-conjecture.

If we consider a fixed $D$, it is seemingly very difficult to find particular $p$ for which 
$K =  \qrD$ is non-$p$-rational. 
Hence the continuing open status of venerated problems like the existence or otherwise of so-called Wall-Sun-Sun or Fibonacci-Wieferich primes, which at least in a na\"{\i}ve sense seem to be inaccessible to current methods of computation~\cite{primegrid}. Indeed, it is this difficulty which lends credibility to the so-called $p$-adic Leopoldt conjecture of Gras~\cite{grascan} and others~\cite{khare}. 

On the other hand, however, if we fix a prime $p$ and search through square-free values of $D$ for which $ \qrD$ is non-$p$-rational, then we easily find huge numbers of such fields. 
Viewed through the latter lens rather than the former, these theorems are not at all surprising.

The primes~$p=2,3$ are excluded only because it would involve otherwise-redundant extra notation and special cases. 
Suffice it to say that similar infinite families may readily be found for these primes using the same approach. 
Moreover, the theorems may be fine-tuned for general $p$, to give vastly larger collections of non-$p$-rational fields by using better bounds on the size of the individual $u_D$ given in propositions \ref{dmn} and \ref{huaest} below, together with equation \eqref{gcdees} and lemmas~\ref{knewp} and~\ref{peeell}.

\subsection{Equivalent formulations for $p$-rationality}
\subsubsection{Notation}
Let $D\geq2$ be a square-free positive integer, and write~$K =  \qrD$ for the real quadratic field with fundamental discriminant $D$ or $4D$ and ring of integers~$\Z_K$. 
Let~$u_K$ be the fundamental unit of~$K$ which is~$>1$ under the embedding of~$K$ into $\R$ sending $\qD$ to a positive number. 
We denote by~$u_D$ the first totally positive power of~$u_K$:
in other words, $u_D=u_K^2$ if and only if~$K$ contains a unit of absolute norm~$-1$; otherwise $u_D = u_K$. 
For any integer~$\ell\geq1$ the~$\ell$-th `dimension' in the tower above~$ \qrD$ is given by $1$ plus the trace of $u_D^\ell$, viz.:
\begin{equation}\label{dees}
d_\ell(D) = u_D^\ell + u_D^{-\ell} + 1. 
\end{equation}
These sequences form a third-order homogeneous (or second-order non-homogeneous) linear recurrence relation and their terms may be charted by a simple variant of the Chebyshev polynomials of the first kind: all of which we outline in appendix~\ref{appx}, along with their generating function. 
Because of how they arose, from now on for convenience we shall refer to the $d_\ell(D)$ as \emph{dimensions} without further ado.

Let~$\cR$ be a Dedekind domain and let~$\cJ \triangleleft \cR$ be any ideal. 
Choose an element~$\epsilon \in \cR$ such that the principal ideal~$(\epsilon) \triangleleft \cR$ is prime to~$\cJ$. 
We denote by~$\ort_\cJ(\epsilon)$ the order of~$\epsilon$ inside the multiplicative group~$\left({\cR} / {\cJ}\right)^\times$. 
The phenomenon with which we are concerned in this note is when the combination~$\{p,K\}$ satisfies the following condition:
\[
\ort_p(u_K) = \ort_{p^2}(u_K) . \tag*{\blah}
\]
More generally, given any integer $r\geq3$ the condition may be modified to be the stronger affirmation that $\ort_p(u_K) = \ort_{p^r}(u_K)$. 
Since the exact value of $r$ is not generally important for our purposes, provided that it is $\geq2$, we shall for the sake of brevity simply refer to them by saying that the combination of $p$ and $K$ satisfies~\blah.

\subsubsection{Properties of the fields $K$ which are equivalent to $p$-rationality}
The finite frequency of the occurrence of such primes as satisfy \blah\ --- in our context for a fixed real quadratic field $K = \qrD$ --- is referred to by B\"ockle \emph{et al} in~\cite{khare} as a `$\bmod\ p$ Leopoldt's conjecture'. 
Indeed, in a series of fascinating papers over many years --- we cite but a few in this paper --- Gras has woven this into a kind of local-global principle for logarithm maps, captured for example in the paper~\cite{grascan}, wherein in \S8 it forms the basis for a conjecture that for a fixed $K$, the condition \blah\ may only be fulfilled by at most a finite number of finite primes $p$. 
Other authors, for example \cite{khare}, weaken this conjecture to applying to a set of density zero. 
In any case the heuristics are far from convincing in either direction at this stage, as the growth of the examples where $p<n$ displays this behaviour for a fixed field $K$ would appear to be on the order of $\log\log n$ \cite[\S8.3]{wash}, and consequently for the same reasons as in~\cite{primegrid} not currently amenable to computer calculations. 

As we noted above, for a fixed real quadratic field~$K$ and a fixed prime $p$ the notion of~$p$-rationality of a real quadratic number field boils down to whether or not~\blah~is satisfied, unless~$p \mid h_K$, the class number of $K$. 
So for all but finitely many primes~$p$ we may say that $K$ is non-$p$-rational if and only if~\blah~is true. 
Since our approach gives us no particular handle on~$h_K$, we shall confine ourselves to calculations which prove~\blah~for particular combinations of~$p$ and~$K$. 

Because of the profusion of terminology for this phenomenon in the literature we list a few alternative ways of stating what `non-$p$-rational' means in this real quadratic field context, together with some references where these relations may be sought. 
This list, sadly, is by no means comprehensive but we refer to Gras' relatively recent article in \cite{grasincram} for most of the context into which these varying categorisations may be placed. 
There is also \cite{lenny}, which while it does not apparently have any direct relationship to our techniques, nor to the others mentioned in the corollary, is nevertheless a curious phenomenon.

\begin{corollary}\label{cuervo} 
Let $p\geq 5$. 
For every real quadratic field $K= \qrD$ and $r$ as in the theorem, the following statements about $p$ are either equivalent to \blah, or are implied by it, as shown in each case.
Consequently in particular they all imply non-$p$-rationality, and are equivalent to it if and only if $p\nmid h_K$. 
\begin{enumerate}
\item ($\iff$\blah) \cite{grascan} the normalised $p$-adic regulator $\mathcal{R}_K$ of $K$ has $p$-valuation $\geq1$

\item ($\iff$\blah) \cite{grascan} the $p$-adic logarithm $\log_p(u_K)$ of $u_K$ is congruent to zero modulo $p^{2}$; 

\item ($\iff$\blah) \cite{wall,sunsun} $p$ is a $k$-Wall-Sun-Sun prime (or generalised Wall-Sun-Sun prime, or Fibonacci-Wieferich prime \cite{Bou1}), where the $k$ refers to the generalised Lucas sequence associated to $K$ being of type $(k,-1)$; 

\item ($\iff$\blah)  \cite{Bou1} the period of the generalised Fibonacci sequence $(U_n(k,-1))$ attached to $K$ is the same modulo $p^2$ as it is modulo $p$; 

\item ($\iff$\blah) \cite[proposition 4.1.1]{ralph} 
for each prime $\fp$ of $K$ lying above $p$, the field completion $K_\fp$ contains a $p$-th root of the image of $u_K$ in $K_\fp$; 

\item ($\iff$\blah) \cite{dutta} (only in the special case where~$K=\Q(\sqrt{2})$): $p$ is a balancing Wieferich prime for $K$;

\item ($\impliedby$\blah) \cite{grascrelle1} the Galois group of the maximal abelian pro-$p$-extension of $K$ unramified outside primes above $p$ contains a cyclic torsion subgroup of order at least $p$; 

\item  ($\impliedby$\blah) \cite[Thm7.14]{wash} the Iwasawa $\nu$-invariant for $K$ is non-zero;

\item  ($\impliedby$\blah)  \cite{NQD} the (Pontrjagin dual of the) Galois cohomology group $H^2(\Gal_{\overline{K}/K},\Z/p\Z)$ is not zero.  \qed
\end{enumerate}
\end{corollary}

\section{Dimension towers above $K$ and the property \blah}
\subsection{Link between dimension towers and \blah}
The link between this theory of $p$-rationality and our dimension towers is the following. 
For the proof, we shall need to invoke a couple of results stated in the next section. 
\begin{proposition}\label{jumpers}
Let~$\ell\geq1$ be the index such that $d_\ell(D)$ contains
the first instance of the prime $p$ as a divisor 
in the tower of dimensions above $ \qrD$, with say $p^r \lVert d_\ell(D)$ for some $r\geq 1$, meaning that $p^r$ is the exact power of $p$ dividing $d_\ell(D)$. 
A sufficient condition for $p,D$ to satisfy \blah\ is that~$p^2 \lvert d_\ell(D)$. 
\end{proposition}
\noindent
This is not a necessary condition: for example, $p=13$ satisfies \blah\ for $D=2$, but no dimension of the form $d_k(2)$ is divisible by $13$. 
Indeed, $13$ belongs to the $\frac{5}{8}$ of all primes which do not appear anywhere as divisors of some $d_k(2)$. 
On the other hand $d_5(2) = 7 \cdot 31^2$ and so $31$ is covered by this lemma, with $\ell=5$ and $r=2$. 

\begin{proof}
Set $p,D,r,\ell$ to be as in the statement of the proposition, so that from \eqref{dees}:
\begin{equation}\label{pKr}
d_\ell(D) = u_D^\ell + u_D^{-\ell} + 1 \equiv 0 \bmod p^r .
\end{equation}
Since $\ell$ is minimal for this property, we know by lemma \ref{legs} that $\ort_p(u_K) = 3\ell$. 
In particular, therefore, $u_D^\ell \not\equiv 1 \bmod p^r$ (if it were then moreover \eqref{pKr} would reduce to $3 \bmod p^r$, a contradiction since $p\neq3$). 
Hence multiplying \eqref{pKr} by $u_D^\ell - 1$ we do not change the $p$-divisibility properties; \emph{a fortiori} since $u_D$ is a global unit, we do not change them by multiplying by $u_D^\ell$ itself either. 
So finally
$$
u_D^{3\ell} - 1 = (u_D^\ell-1) u_D^\ell (u_D^\ell + u_D^{-\ell} + 1) \equiv 0 \bmod p^r .
$$
But by \eqref{ordiv}, $\gcd(p,3\ell) = 1$, so by the definition and basic properties of the $p$-adic logarithm \cite[\S5]{wash}:
$$
\log_p(u_D) = 3^{-1}\ell^{-1} \log_p(u_D^{3\ell} ) = \sum_{n\geq1} \frac{(-1)^n}{n} ( u_D^{3\ell} - 1 )^n \equiv 0 \bmod p^r .
$$
But $\log_p(u_K) = \log_p(u_D)$ or $2^{-1}\log_p(u_D)$ and so since $p\neq2$ we are done by the second equivalence in corollary~\ref{cuervo}. 
\end{proof}

\subsection{Basic lemmas on `dimension towers' above $K$}\label{scoo}
Next we show how the dimensions fit together inside the partial order above each fixed real quadratic field $ \qrD$. 
The basic set-up is taken from \cite[\S3]{AFMY} and \cite[appendix A]{BGM}. 
These papers are about `crystallography' in complex $d$-dimensional Hilbert space, which in principle does not concern us here; but nevertheless we retain the terminology so as to be able to access the partial order, by divisibility, of positive integers (what we continue to refer to simply as `dimensions') $d_\ell(D) \geq 4$ attached to any chosen fixed real quadratic field $K= \qrD$, as we now describe.

Recall that given any integer~$n\geq4$ we defined 
$\D(n)$ to be the square-free part of $ (n+1)(n-3) $. 
Up to a canonical factor of~4, $\D$ maps the set of dimensions~$\N \setminus \{0,1,2, 3\}$ 
onto the set of fundamental discriminants of all proper real quadratic fields, since
for every dimension $d = d_\ell(D)\geq4$ in~\eqref{dees}, 
\[
(d-1)^2 - f^2 D = 4,
\]
the Pellian equation defining the integer $f\geq1$. 
The solution $(d-1) + f\qD$ must be some positive power of the fundamental unit $u_K$, since the $\Z$-rank of the unit group of $\Z_K$ is $1$. 
See the table in the appendix to this paper in order to form some intuition for how this plays out in practice. 
Indeed, every power of $u_D$ gives such a unique dimension via the relation \eqref{dees}, and so as observed in \cite[\S3]{AFMY} 
it follows that above any real quadratic field $K$ we 
obtain an infinite tower of dimensions $\{ d_\ell(D) \}_{\ell\geq1}$ which is totally ordered and monotonically increasing in ordinary real absolute value by the index~$\ell$, and also partially ordered by divisibility. 
For, as shown in \cite[proposition A.5]{BGM} for every pair~$\ell,{\ell'}\in\N$ for which~$\frac{\ell}{{\ell'}}$ in its lowest 
terms has no net power of $3$ dividing it, 
\begin{equation}\label{gcdees}
d_{\gcd(\ell,{\ell'})}(D) = \gcd(d_\ell(D),d_{\ell'}(D)) . 
\end{equation}
It is this \emph{strong divisibility} structure --- modulated by the shears occurring at each power of $3$ --- which enables the observations of this note. 
For completeness we mention that when the power of~$3$ dividing $\ell$ is different from that dividing~${\ell'}$, then $\gcd(d_\ell,d_{\ell'}) \in \{ 1, 3 \}$: this messy complication led us to avoid treating the case of $p=3$ in this paper. 

We briefly recall some other results. 
Write $\leg{a}{p}$ for the usual Legendre symbol mod $p\neq2$. 
For simplicity we continue to assume that $\gcd(p,6D)=1$. 
First, by lemma 3 and the remark following the theorem of \cite{tbang}, using the usual connection between Chebyshev polynomials of the first kind, the hyperbolic cosine and powers of $u_D$, 
\begin{equation}\label{ordiv}
\ort_p(u_D)\ \Big\lvert  \ p-\leg{D}{p} .
\end{equation}
\noindent
The next lemma is explained in the text in section 3 of \cite{AFMY}. 
\begin{lemma}\label{threase}
Fix $D$ and the tower of fields $\{ d_k(D) \}_{k\geq1}$ as above. 
There are infinitely many `co-prime' sub-towers $\{ d_{3^sk}(D) \}_{k\geq1; 3\nmid k}$ --- in the sense that with the possible exception of the prime $3$, the prime divisors of one are disjoint from those of any other --- whose base points are indexed by the non-negative integer powers~$3^s$. \qed
\end{lemma}
\noindent
In other words, there is an infinite series of sub-sequences indexed by $s\in\N\cup\{0\}$, which we could also write as $\{ d_{3^sk} \}_{k \in \N \setminus 3\N}$, which may all effectively be re-based as though $3^s = 1$ for the purposes of this paper and analysed separately from one another. 
See equations (113) and (114) in \cite{BGM}, together with the appendix to this paper. 

\begin{lemma}[\cite{AFMY,kopp}]\label{legs}
Let $p\geq5$ be a prime, and let $D\geq2$ be square-free, co-prime to $p$. 
Then $p$ appears in the dimension tower above $ \qrD$ --- that is to say, there is an $\ell\geq1$ and a power $r\geq1$ for which $p^r\ \lVert\ d_\ell(D)$ --- if and only if $\ort_p(u_D)$ is divisible by 3. 
In this case, $p \equiv \leg{D}{p} \bmod 3$ and $\leg{D}{p} = \leg{p}{3}$. 
Indeed for the minimum such $\ell$, $\ort_p(u_D) = \ort_{p^r}(u_D) = 3\ell$. \qed
\end{lemma}

\noindent 
For the next lemma, for technical reasons explained in the proof of the proposition cited from \cite{BGM}, we must exclude cases where $d_1(D)$ is of the form $2^N+2$ for some $N\geq1$ and where $\ell=2$, so the problem is with $d_2(D) = 2^{N+1}(2^{N-1}+1)$: see the formulae in \S3 of \cite{AFMY}. 
For example, we would need this extra stipulation in order to exclude pathological cases like $d_1(5) = 4$ and $d_2(5) = 8$; or indeed $d_1(77) = 10$ and $d_2(77) = 80$, wherein in ascending one rung of the ladder of dimensions we do not in fact introduce any new primes. 
However in our context we lose nothing by this, since we know by the most elementary case of Catalan's conjecture --- now Mihaelescu's theorem --- that no proper prime power other than $9$ can be of the form $2^M+1$; hence numbers of the form $d_2(D)$ as posited do not occur in the sequences of our theorems. 

\begin{lemma}[Proposition A.2 of \cite{BGM}]\label{knewp}
Excluding the case just mentioned, for every $\ell \geq 1$ and every $D$ the dimension $d_\ell(D)$ introduces at least one `new' prime number~$q$ into the tower above $ \qrD$. 
That is, $q$ divides into $d_\ell(D)$ but not into $d_k(D)$ for any $k < \ell$. \qed
\end{lemma}

\begin{lemma}[Proposition A.3 of \cite{BGM}]\label{peeell}
Let $p\geq5$ be a prime and let $D\geq2$, not divisible by $p$, be square-free. 
Assume that $p,D$ satisfy the conditions of lemma~\ref{legs}, and let $k,r \geq 1$ be such that $p^r \lVert d_k(D)$. 
Then~$p^{r+1} \lVert d_{pk}(D)$. 
Furthermore $pk$ is the lowest such integer. \qed
\end{lemma}

\section{Proofs of the theorems}

\subsection{Proof of theorem \ref{peequeue}}\label{proof1}
We invoke lemma \ref{threase} at the outset, in order that the analysis may be self-contained. 
So without loss of generality we assume that we have accounted for any powers of $3$ which divide the index $\ell$ of $d=p^rq^t = d_\ell(D)$ within its sub-tower above $ \qrD$. 
They are thereby assumed to be factored out so that the root of the sub-tower is just called $d_1(D)$. 
What follows remains an entirely general argument under such a hypothesis, thanks to eqn.~\eqref{gcdees}. 

By proposition \ref{jumpers}, we need to show that when $p,q,r,t,D$ satisfy the conditions in the theorem, they force $d=p^rq^t$ to be of the form $d_\ell(D)$ for some $\ell$ which is not divisible by $p$; or if it is, $r\geq3$ so that the combination $p,K$ still satisfies \blah. 
To this end, observe by lemma \ref{threase} that if we set $D = \D(p^r q^t)$, and then work out for which value of $k\in\N$ we have $d = d_k(D)$, then $k$ must either be a power of $3$, which by the observations made above we re-base to 1; or else $k$ must itself be prime. 
Otherwise, given that $d_1(D)\geq4$, too many new primes would have had to have been introduced in going up the tower of divisibility relations \eqref{gcdees}, by lemma \ref{knewp}. 

If $k=1$ we have the answer we need, by proposition \ref{jumpers}; on the other hand if $q^t=1$ then we once again are done by proposition \ref{jumpers} and lemma \ref{knewp} since it then must be the case that $d_1(D) = p^r$. 
So we may assume either $k=p$, $k=q$, or $k$ is some other prime, and moreover that $q^t \geq 2$. 

If $k=p$ then by lemma \ref{peeell}, $d_1(D) = p^{r-1}$ satisfies \blah\ for $r-1$ rather than $r$ (note that $q$ is the `new' prime for level $p$, so it cannot divide into $d_1(D)$ here). 
Provided $r\geq3$ this is not a problem, as it still satisfies a non-$p$-rationality criterion. 
So without loss of generality we may assume that $r=2$; then the sequence of dimensions for this particular $D$ would have to begin with $d_1(D) = p$ and $d_p(D) = p^2q^t$. 
This is perfectly possible. 
For example, if $D=2$ and we choose $p=d_1(2) = 7$ then $d_7(2) = 7^2 \cdot 4663$, so setting $q$ to be the prime $4663$ with $t=1$, we have constructed such a case. 
So we must bound $q$ somehow to preclude this. 
An easy bound, which is too loose to be useful for anything other than proving this theorem, comes from proposition \ref{dmn} below, where it is shown that for all $D$ (the case for $D=5$ is true with a slightly modified lower bound but it is not relevant here), for all $\ell\geq1$ and for all $n\geq2$:
\begin{equation}\label{toodle}
\frac{2^n}{3^n} < \frac{d_{\ell n}(D)}{d_\ell(D)^n} < 1.
\end{equation}
The idea here is to use $\ell=1$ and $n=p$ and $d_\ell(D)=p$ also; so in particular, if we stipulate that 
\[
t < p\log_q(\frac{2}{3}) + (p-2)\log_q p
\]
where this time (and only here) the notation $\log_q$ simply means the ordinary logarithm to the base $q$, then it follows immediately that $q^t < \frac{2^p}{3^p}p^{p-2}$ and so if $d_1(D) = p$ and $d_p(D) = p^2q^t$ as we have assumed, we get 
\[
d_p(D)  = p^2q^t < \frac{2^p}{3^p}p^p =  \frac{2^p}{3^p} d_1(D)^p  ,
\]
a contradiction to the lower bound in \eqref{toodle} with $\ell=1$ and $n=p$. 
Indeed, in the example given above of $d_7(2)$, our bound would force $t<0.815787..$ and so we would already have excluded this case, as desired. 

If $k=q$, then $d_1(D) = q^{t-1}$ again by lemma \ref{peeell} and the following argument: the other possibilities would be $d_1(D) = p^rq^{t-1}$, which would then violate lemma \ref{knewp} at the $q$-th level; or else $d_1(D) = p^r$ which would violate \eqref{ordiv} for the first appearance of some power of~$q$. 
Hence by proposition \ref{jumpers} once more we have satisfied the requirements. 
Note that if $q=3$ we are forced to rebase as above, so we may effectively exclude $q=3$ here. 

So finally we are left with the case where $k$ is some prime other than $3,p,q$. 
For example $D=2$, $p=31$, $q^t=7$ so $d_1(2) = 7$ and for $k=5$ we get $d_5(2) = 7 \cdot 31^2$: this is an instance of the phenomenon we are seeking. 
This is straightforward if $d_1(D) = q^t$ or $d_1(D) = p^r$, for then we clearly shall have satisfied the requirements for applying the result of proposition \ref{jumpers}. 
But no other combinations can happen, by lemmas \ref{knewp} and \ref{peeell}, since $k$ is assumed co-prime to $pq$ and since $d = d_k(D) = p^rq^t$ only contains two distinct primes. 

This completes the proof of theorem~\ref{peequeue}. 
\qed

\vskip 1mm

The statement in the title and abstract --- that the set of distinct such $D$ thus generated can be made to be infinite --- follows immediately from the rest, if for example we just take the sequence of dimensions $p^2,p^3,p^4,p^5,\ldots$. 
But consider two arbitrary elements of a sequence $p^rq^t$, where $p$ and $q$ are kept fixed and $r$ and $t$ allowed to grow. 
Suppose that $\D(p^rq^t) = \D(p^{r'}q^{t'}) = D$, say, for some $r<r'$ and $t<t'$. 
Let $p^rq^t = d_k(D)$ and $p^{r'}q^{t'} = d_{k'}(D)$ define indices $k,k'$. 
By \eqref{gcdees} we know that $k \mid k'$; by lemma \ref{threase} since $p\neq q$ they must be in the same sub-tower (ie with base point $d_{3^s}(D)$ corresponding to the same power $3^s$) and so it follows from lemma \ref{knewp} that a new prime must be introduced in going from $d_k(D)$ up to $d_{k'}(D)$, and this is clearly not the case. 
So any such collection of fields $ \qrD$ based on the same two primes $p,q$, with a non-decreasing sequence of $r$ and $t$ respectively, will give rise to an infinite collection of fields. 

As an example, let $[ t_j : j \geq 3 ]$ be an infinite sequence of non-negative integers, with or without repeated values. 
Now take any prime $p\geq 5$ and any prime $q\neq p$: the set of dimensions $\{ d^{(r)} := p^r q^{t_r} : r \geq 3 \}$, for the reasons given in the proof of the theorem, is such a set giving rise in a 1-1 mapping to a set of square-free $D^{(r)} := \D( d^{(r)} )$, all of which are distinct and all of which give distinct non-$p$-rational real quadratic fields $\Q(D^{(r)})$. 

\subsection{Proof of theorem \ref{pqr}}\label{proof2}
Just for orientation's sake, the initial case $n=1$ may be shown by taking the case $p^2$ of theorem \ref{peequeue}, with $t=0$. 
So we consider a general $n\geq2$ now, and follow a similar argument to that of the previous proof. 
However this time it is easier because the restriction upon $t$ in the theorem \ref{peequeue} is replaced by a much stronger restriction on the overall shape of the numbers $d^{(\bfr)} = \scriptstyle{ \prod\limits_{j=1}^n p_j^{r_j} } $,  
namely that all of the individual powers $r_j$ be $\geq n+1$. 
This makes the use of lemmas \ref{knewp} and \ref{peeell} far more straightforward. 

Suppose the statement is false for some fixed set of primes $p_1,\ldots,p_n$. 
Let $\bfr = (r_j)_{1\leq j \leq n}$ be a specific choice of index vector for which the field $K^{(\bfr)}$ defined in the statement of the theorem fails to satisfy \blah\ for some particular $c$, $1\leq c\leq n$. 
Recall that we apply the `discriminant' map $\D$ to the dimension defined as $d^{(\bfr)} = \scriptstyle{ \prod\limits_{j=1}^n p_j^{r_j} }$, we obtain some $\D(d^{(\bfr)}) = D$ and finally the quadratic field of interest which is $K^{(\bfr)} = \qrD$. 
Let $\kappa$ denote the value of $\ell$ corresponding to the dimension $d^{(\bfr)}$ itself: that is to say, 
$d_\kappa(D) = d^{(\bfr)}$. 

\def\vp{v_{p_c}}

Since the argument becomes easier the larger any $r_c$ is, we lose nothing by assuming that $r_c = n+1$. 
Let $k_c \in \N$ be the minimal value of $\ell$ for which $d_\ell(D)$ is divisible by some power of $p_c$. 
Let~$\vp \colon \Q \lra \Z$ denote the normalised $p_c$-adic valuation: in particular then, by proposition \ref{jumpers} and our assumption that \blah\ is not satisfied for $p_c$, $\vp( d_{k_c}(D) ) = 1$. 
By lemma \ref{peeell} this forces the following structure among the dimensions in the component of the $p_c$-adic tower leading up to $d^{(\bfr)}$:
\begin{equation}\label{hada}
\scriptstyle
\vp(d_{k_1}(D)) = 1, \ \vp(d_{p_ck_1}(D)) = 2, \ \vp(d_{p_c^2k_1}(D)) = 3, \ \ldots, \ \vp(d_{p_c^{n-1}k_1}(D)) = n, \ \vp(d_{p_c^{n}k_1}(D)) = n+1 .
\end{equation}
Now by lemma \ref{peeell}, ${p_c^{n+1}k_1}$ must divide into $\kappa$; 
hence since we have excluded $3$ from the $p_j$, by \eqref{gcdees} it follows that $d_{p_c^{n}k_1}(D) \left\lvert\right. d^{(\bfr)}$. 
On the other hand we also excluded $2$ from the $p_j$ so that by lemma \ref{knewp}, there must be at least one `new' prime $q_i$ coming in to each of these successive levels in \eqref{hada}, for $2\leq i \leq n+1$. 
But then we would have at least $n$ new mutually distinct prime divisors $\{ q_i \}_i$ on the way from $d_{k_1}(D)$ up to $d_{p_c^{n}k_1}(D)$ and so \emph{a fortiori} also dividing into~$d^{(\bfr)}$. 
This is a contradiction, since $d^{(\bfr)}$ was defined to have only $n$ distinct prime factors.

This completes the proof of theorem~\ref{pqr}, other than the statement regarding the collection having no repeated elements (and therefore being infinite), which follows immediately via the same argument as we explained above for the infinite cardinality of the sets generated from theorem \ref{peequeue}. 
\qed

\subsubsection*{Acknowledgements}
As mentioned above, I would like to express my gratitude to John Coates, who sadly passed away in the time since I began this work, for his mentoring throughout my life. 
I also owe a great debt to Georges Gras, whose detailed comments on an earlier draft of this paper were indispensable. 
Thanks to Gene Kopp, Jeff Lagarias and Chris Smyth for helping me to set this research on a firm footing in the early days. 

In addition I would like to thank John Cannon and his team for MAGMA~\cite{Magma}, which has been central to everything done here. 
Moreover, without the prompting of Terry Rudolph and then the generous collaborations over many years with Marcus Appleby, Ingemar Bengtsson, Markus Grassl and Mike Harrison, I never would have found the problem with which this paper is concerned. 
Finally I would like to acknowledge gratefully the ongoing generosity of Myungshik Kim and the QOLS group at Imperial College London.

\appendix
\section{Modified Chebyshev polynomials and dimension towers}
\label{appx}

In order that this presentation at least be somewhat self-contained, we very briefly cover some of the techniques upon which we have drawn in the main text, from references \cite[\S3]{AFMY} and \cite[Appendix A]{BGM}. 

\subsection{The dimensions in simple terms}
We begin by explaining the `dimensions' to which we are constantly referring, in a down-to-earth form. 

It is an obvious yet powerful fact that with the exception\footnote{These come in only in the negative norm cases~$D=2: u_K = 1+\sqrt{2}$ and~$D=5: u_K = \frac{1+\sqrt{5}}{2}$. As soon as we introduce \emph{all} powers of these into the equations the correspondence will be 2-1 --- rather than 1-1 --- other than for these two initial points.} of the initial two points~$\{ \frac{1}{2} , 1 \}$, we may put the set~$\frac{1}{2}\N$ in one-to-one correspondence with the rational components of (unique) powers of~$u_D$ for appropriate~$K =  \qrD$, for example by using the formalism of our dimensions $ \frac{1}{2}\Tr u_D^\ell \longleftrightarrow u_D^\ell \longleftrightarrow \Tr u_D^\ell +1 =:  d_\ell(D)  \longleftrightarrow \{\ell,D\}$. 
For example:
\begin{eqnarray*}
\frac{3}{2} \mapsto u_5 = u_{\Q(\sqrt{5})}^2 = \frac{3+\sqrt{5}}{2} \longleftrightarrow d_1(5) = 2\times\frac{3}{2}+1 &=& 4 \\
\frac{4}{2} = 2 \mapsto u_3 = u_{\Q(\sqrt{3})} = 2 + \sqrt{3} \longleftrightarrow d_1(3) = 2\times\frac{4}{2}+1 &=& 5 \\
\frac{5}{2} \mapsto u_{21} = u_{\Q(\sqrt{21})} = \frac{5+\sqrt{21}}{2} \longleftrightarrow d_1(21) = 2\times\frac{5}{2} + 1 &=& 6 \\
\frac{6}{2} = 3 \mapsto u_2 = u_{\Q(\sqrt{2})}^2 = 3 + 2\sqrt{2} \longleftrightarrow d_1(2) = 2\times\frac{6}{2}+1 &=& 7 \\
\frac{7}{2} \mapsto u_5^2 = u_{\Q(\sqrt{5})}^4 = \frac{7 + 3\sqrt{5}}{2} \longleftrightarrow d_2(5) = 2\times\frac{7}{2}+1 &=& 8 \\
\frac{8}{2} = 4 \mapsto u_{15} = u_{\Q(\sqrt{15})} = 4 + \sqrt{15} \longleftrightarrow d_1(15) = 2\times\frac{8}{2}+1 &=& 9 \\
\frac{9}{2} \mapsto u_{77} = u_{\Q(\sqrt{77})} = \frac{9 + \sqrt{77}}{2} \longleftrightarrow d_1(77) = 2\times\frac{9}{2}+1 &=& 10 \\
\frac{10}{2} = 5 \mapsto u_6 = u_{\Q(\sqrt{6})} = 5 + 2\sqrt{6} \longleftrightarrow d_1(6) = 2\times\frac{10}{2}+1 &=& 11\ \ldots \hbox {\ et\ cetera} 
\end{eqnarray*}

\subsection{Chebyshev polynomials and navigating the towers above a fixed $K$}
The final short sub-section is taken almost verbatim from \S A.2 of \cite{BGM}. 
In order to navigate among the dimensions in the same `tower complex' 
above a fixed~$K$, 
we introduced in~\cite[\S3]{AFMY} a family of 
functions which enable us to calculate any 
dimension~$d_{n\ell}(D)$ given the value of~$d_\ell(D)$, 
akin to the role played by the Chebyshev polynomials
$T_n(x), U_m(x)$ for the trigonometric functions. 
These modified Chebyshev polynomials of the first kind are defined by: 
\begin{equation}\label{shifty}
\chsh{n}{x} = 1+2\cheb{n}{\tfrac{x-1}{2}},
\end{equation}
where~$T_n(x)$ is the usual Chebyshev 
polynomial of the first kind of degree~$n$. 
Like the ordinary Chebyshev polynomials, the~$\chsh{n}{x}$ also satisfy the recursive relation
\begin{equation}\label{nest}
\chsh{m}{\chsh{n}{x}} = \chsh{mn}{x} \hbox{\ for\ every\ }m,n\geq0.
\end{equation}
The~$\chshh{n}$ are independent of~$D$ so once we 
know~$d_0(D) = 3$ (which is true trivially, by applying the definitions formally, for all~$D$) 
and~$d_1(D)$ --- which requires that we know the 
fundamental unit of~$K$ --- we know 
all~$d_\ell(D)=\chsh{\ell}{d_1}$. 
The first few polynomials are~$\chsh{0}{x}=3$, $\chsh{1}{x}=x$, 
$\chsh{2}{x}=x^2-2x$, $\chsh{3}{x}=x^3-3x^2+3$. 
The defining recursion $\cheb{n}{x} = 2x\cheb{n-1}{x} - \cheb{n-2}{x}$ 
for the~$\chebh{n}$ yields for the~$\chshh{n}$:
\begin{equation}\label{shiftycur}
\chsh{n}{x} = x\chsh{n-1}{x} - x\chsh{n-2}{x} + \chsh{n-3}{x}.
\end{equation}
We recall lemma \ref{threase}, where we may re-base the tower to start from any dimension $d_\ell(D)$. 
Then \eqref{shiftycur} gives recurrence relations just as if we had started with $d_1(D)$ instead, viz.: 
$$
d_{n\ell} = d_\ell ( d_{(n-1)\ell} - d_{(n-2)\ell} ) + d_{(n-3)\ell} ; 
$$
which in particular when~$\ell=1$ becomes: 
$$
d_n = d_1 \cdot (d_{n-1} - d_{n-2}) + d_{n-3} . 
$$
This is a third-order homogeneous representation of the recursion; on the other hand we may reduce it to a non-homogeneous second order recursion by rearranging thus: 
for any~$n\geq2$ (and noting once again that~$d_0(D) = \Tr u_D^0 + 1 = 3$ for any~$D$):
$$
d_n = (d_1-1) d_{n-1} - d_{n-2} - (d_1-3) . 
$$
Note this is also true of the~$d_{n\ell}$ expressions above, with appropriate amendments to the indices for~$d_\ell$ replacing~$d_1$ and yielding thereby:
$$
d_{n\ell} = ( d_\ell - 1) \cdot d_{(n-1)\ell} - d_{(n-2)\ell} - (d_\ell - 3 ) .
$$ 

\subsubsection{Generating function for the $d_\ell(D)$ as a logarithmic derivative}
Purely for interest's sake, in the hope it will give extra handles on the material for some readers, we mention the generating function for the $d_\ell(D)$. 

\begin{proposition}\label{secordlin}
Fix~$D$. 
Denoting the minimal polynomial of~$u_D$ by 
$$
m_1(x) = (x - u_D)(x-u_D^{-1})  =  x^2 - (d_1 - 1)x + 1 ,
$$
the (ordinary) generating function~$\Gamma(x)$ for the~$d_k$ may be constructed using standard techniques, to be: 
$$
\Gamma(x)  =    \frac{ 1 }{  1  -  x }    +    \frac{ u_D }{  u_D  -  x }  +  \frac{ u_D^{-1} }{  u_D^{-1}   -   x }    =   \frac{ 3 - 2 d_1 x + d_1 x^2 }{ 1 - d_1 x + d_1 x^2 - x^3 }     \qed.
$$
\end{proposition}

We note that this generating function~$\Gamma(x)$ has the form of a logarithmic derivative of sorts, in that defining a function~$F(x) = (1-x) m_1(x)$ : 
$$
\frac{1}{x}\Gamma\left({{\frac{1}{x}}}\right) = \frac{F'(x)}{F(x)} , 
$$
Write $t_k=d_k-1$, the trace of $u_D$ itself.  
Then this and the generating function for the traces, are the only two series of the form~$t_k+\lambda$ for which this functional equation is true, for any~$\lambda\in\R$: for it is merely a formal consequence of the $k$-th term being defined as the sum of the~$k$-th powers of the roots of a polynomial.   

\subsubsection{Auxiliary results on comparative sizes of dimensions in the towers}
The following results are taken from \cite[Appendix A]{BGM}. 
We have used them in the proof of theorem \ref{peequeue}. 
\begin{proposition}\cite[Proposition A.4]{BGM}\label{dmn}
Let~$D\ne5$ be any square-free positive integer\footnote{There is an easy minor modification for~$D=5$ for both of these results (which is only needed 
for the first few dimensions anyway) but we do not go into this here.}. 
For all $\ell\geq1$ and~$n\geq2$:
\[
d_\ell^n > d_{n \ell} > \left(\frac{2}{3}\right)^n d_\ell^n .
\qed\]
\end{proposition}
\noindent
This relies on the following result, which essentially goes back to \cite{hua}. 

\def\D{\mathcal{D}}
\begin{proposition}\label{huaest}\cite[Lemma A.6]{BGM}
Let $D>1, D\neq5$ be a square-free integer and let $\Delta_D := f^2
D$, for $f \in \{1,2\}$, be the (`fundamental') discriminant of the
real quadratic field $ \qrD$. That is, $f=1$ if $D\equiv1\bmod4$
and $f=2$ if $D\equiv2$ or $3\bmod4$. Write $\D$ for the positive
real number $\sqrt{\Delta_D} = f\qD$. Let $e$ be the base of the
natural logarithm. Then with the usual notation for the integer part
of a real number,
$$
\lfloor \D \rfloor < \frac{\D+\sqrt{\D^2-4}}{2} \leq u_K < ( e \D )^\D . \qed
$$
\end{proposition}

\end{document}